\documentclass[12pt,a4paper]{article}
\usepackage{amssymb,amsmath,palatino}
\usepackage[mathcal]{euler}

\def\ket#1{\left| #1\right\rangle}

\def\nn{\nonumber}

\newcommand{\IR}{{\mathbb R}}

\DeclareMathOperator{\Mat}{Mat}
\DeclareMathOperator{\Tr}{Tr}
\newcommand{\C}{\mathcal{C}}
\newcommand{\M}{\mathcal{M}}
\newcommand{\R}{\mathbb{R}}
\newcommand{\cs}{\mbox{\upshape C}\ensuremath{{}^*}}
\newcommand{\K}{\mathcal K}
\newcommand{\T}{\mathbb{T}}
\newcommand{\Z}{\mathbb Z}
\newcommand{\into}{\hookrightarrow}
\newcommand{\co}{\mathbb{C}}
\newcommand{\braket}[3]{\left<#1\vphantom{#2}\right|#2\left|#3\vphantom{#2}\right>}

\newcommand{\norm}[1]{\lVert#1\rVert}

\newcommand{\abs}[1]{\lvert#1\rvert}

\DeclareMathOperator{\Spec}{Spec}
\newcommand{\Li}{\mathcal L}
\newcommand{\Hi}{\mathcal H}
\DeclareMathOperator{\ch}{ch}
\DeclareMathOperator{\tr}{tr}
\newcommand{\xz}{\widetilde x_0}
\DeclareMathOperator{\SU}{SU}
\DeclareMathOperator{\SO}{SO}
\DeclareMathOperator{\Aut}{Aut}
\newcommand{\addots}{\mathinner{\mkern-1mu
\raise0.5pt\vbox{\kern7pt\hbox{.}}\mkern2mu\raise4pt\hbox{.}\mkern2mu
\raise7.5pt\hbox{.}\mkern-1mu}}
\newcommand{\Szymanski}{Szyma{\'n}ski}
\newtheorem{proposition}{Proposition}[section]  %proposition
\newcommand{\bprop}{\medskip\begin{proposition} ~~\\ \it}
\newcommand{\eprop}{\end{proposition} \hfill \medskip \\}
\numberwithin{equation}{section}

\begin{document}

\title{Fredholm Modules\\ for
       Quantum Euclidean Spheres}
\date{October 2002}

\author{Eli Hawkins$\strut^{1}$, \ Giovanni Landi$\strut^{2}$ \\[15pt]
$\strut^{1}$Scuola Internazionale Superiore di Studi Avanzati \\
Via Beirut 4, I-34014 Trieste\\
\texttt{\normalsize mrmuon@mac.com} \\[5pt]
$\strut^{2}$Dipartimento di Scienze Matematiche, Universit\`a di
Trieste\\
Via Valerio 12/b, I-34127 Trieste \\
and INFN, Sezione di Napoli, Napoli\\
\texttt{\normalsize landi@univ.trieste.it}}

\maketitle
\begin{abstract}
The quantum Euclidean spheres, $S_q^{N-1}$, are (noncommutative)
homogeneous spaces of quantum orthogonal groups, $\SO_q(N)$.
The $*$-algebra $A(S^{N-1}_q)$ of polynomial functions on each of
these is given by generators and relations which can be expressed in 
terms of a self-adjoint, unipotent matrix.
We explicitly construct complete sets of generators for the $K$-theory (by
nontrivial self-adjoint idempotents and unitaries) and 
the $K$-homology (by nontrivial Fredholm modules) of the spheres
$S_q^{N-1}$. We also construct the corresponding Chern characters in cyclic
homology and cohomology and compute the  pairing of $K$-theory with
$K$-hom\-o\-logy. On odd spheres (i.~e., for $N$ even) we exhibit unbounded
Fredholm modules by means of a natural unbounded operator $D$
which, while failing to have compact resolvent, has bounded 
commutators with all elements in the algebra $A(S^{N-1}_q)$.
\end{abstract}

\vfill
\textsf{\small math.KT/0210139, SISSA/FM 66/2002, DSM-QM/528, ESI Vienna 1220 (2002)}\\[1ex]
\textit{\small 2000 Mathematics Subject Classification. 19D55; Secondary 20G42, 58B34.}

\thispagestyle{empty}
\newpage

\section{Introduction}
There exists a growing literature devoted to the study of examples of
``quantum" and ``noncommutative" spaces. In this paper we shall
dissect one
class of  these, the so called quantum Euclidean spheres $S^{N-1}_q$.
They were first introduced in \cite{FadResTak89} as homogeneous spaces
of quantum orthogonal groups $\SO_q(N)$, which are $R$-matrix
deformations of
the usual orthogonal groups $\SO(N)$.

We shall regard the spheres $S^{N-1}_q$ as
``noncommutative real affine varieties"\/. For such an object, $X$,
the
algebra $A(X)$ is a finitely presented $*$-algebra. In contrast with
classical algebraic geometry, there does not in general exist a
topological point set $X$. Nevertheless, we regard $X$ as a
noncommutative space and $A(X)$ as the algebra of polynomial functions
on $X$.

In the classical case, one can consider the algebra of continuous
functions on the underlying topological space of an affine variety. If
$X$ is bounded, then this is a \cs-algebra and is the completion of
$A(X)$. In general, one defines $\C(X)$ to be the \cs-algebraic
completion of the $*$-algebra $A(X)$. To construct this, one
considers all
possible $*$-representations of $A(X)$ on a countably
infinite-dimensional Hilbert space. Then the norm on $A(X)$ is
defined as the
supremum of the norms in all these representations and the
\cs-algebra $\C(X)$ is
the
completion of $A(X)$ with respect to this norm. The \cs-algebra
$\C(X)$ has the
universal property that any \hbox{$*$-homomorphism} from $A(X)$ to a
separable
\cs-algebra factors through $\C(X)$. In particular, any
\hbox{$*$-representation} of $A(X)$ extends to a representation of $\C(X)$.

For the noncommutative spaces at hand, the algebra $A(S^{N-1}_q)$
was described in \cite{LanMad01} by means of a suitable self-adjoint
idempotent
(a matrix of
functions whose square is itself).
In the present paper we improve on this by giving a clearer and nicer
presentation of $A(S^{N-1}_q)$ in terms of a self-adjoint unipotent (a
matrix of functions whose square is the identity) which is defined
recursively. We then exhibit
all representations of the algebra $A(S^{N-1}_q)$ which in turn
extend to
the \cs-algebra $\C(S^{N-1}_q)$.

The core of the paper is the study of generators of the  $K$-homology  and
$K$-theory of the spheres $S^{N-1}_q$.  The $K$-theory classes will be given by
means of self-adjoint idempotents (naturally associated with  the 
aforementioned
unipotents) and of unitaries in algebras of matrices over $A(S^{N-1}_q)$.
The $K$-homology classes will  be given as (homotopy classes of) suitable
$1$-summable Fredholm modules.
\\
For odd spheres (i.e. for $N$ even) the odd $K$-homology generators are first given in terms of unbounded Fredholm
modules. These are given by means of a natural unbounded operator $D$ which,
while failing to have compact resolvent, has bounded commutators with all
elements in the algebra $A(S^{2n+1}_q)$.

On the way to computing the  pairing of $K$-theory with
$K$-homology, we will exhibit the Chern characters of the generators
of the
$K$-theory and $K$-homology,
as elements in the cyclic homology $HC_*[A(S^{N-1}_q)]$ and cyclic
cohomology $HC^*[A(S^{N-1}_q)]$ respectively.

Needless to say, the pairing is integral (it comes from a
noncommutative index
theorem). Furthermore, the non-vanishing of the pairing will testify
to the
non-triviality  of the elements that we construct in both
$K$-homology $K$-theory.

\section{Quantum Euclidean Spheres}
\label{se:qes}
As we have already mentioned, the quantum Euclidean spheres,
$S^{N-1}_q$, were
introduced in \cite{FadResTak89} as
quantum homogeneous spaces of the quantum orthogonal groups,
$\SO_q(N)$, the latter
being
$R$-matrix deformations of usual orthogonal groups $\SO(N)$.

We shall briefly recall
that construction \cite{Ogi92,CerFioMad00,LanMad01}.
Let us start with the quantum
Euclidean space $\IR^N_q$. The $*$-algebra $A(\IR^N_q)$ of its
polynomial
``functions''
is generated by elements $\{\xz=\xz^*, x_i, x_i^*,  i=1,\dots,n \}$
for
$N=2n+1$ while for $N=2n$ there is no $\xz$ (the reason for the
notation $\xz$ will be
clear shortly). These generators fulfill  commutation relations,
\begin{subequations}
\label{old.core}
\begin{align}
\begin{split}
x_i x_j &= q x_j x_i ,  \quad 0\leq i<j \leq n, \\
x_i^* x_j &= q x_j x_i^*, \quad
i\neq j ,
\end{split}
\label{old.core1} \\
   [x_i, x_i^*] &=
\begin{cases}
[ (1 - q^{-2}) /  (1 + q^{-2\rho_{i-1}} )] r^2_{i-1} & i > 1
,\\
0 &  i=1, N=2n
,\\
(1-q^{-1})  \xz^2 &  i=1, N=2n+1
,\end{cases} \label{old.core2}
\end{align}
\end{subequations}
with $\rho_{i}=\frac{1}{2}-i$ or $\rho_{i}=1 -i, \; i=1, \dots, n$,
according to
whether
$N$ is odd or even
respectively. The elements $r^2_i, \, i=1, \dots, n$, are given by,
\begin{equation}
r^2_i=q^{-2\rho_i} x^i (x^i)^* + q^{-2\rho_{i-1}} x^{i-1} (x^{i-1})^*
\dots + (x^{i-1})^* x^{i-1} + (x^i)^* x^i ,
\label{defrad}
\end{equation}
and one can prove that $r^2 \equiv r_n^2 $ is central.
\\
In the
classical ($q=1$) case, these
$r_i^2$'s  are simply sums of
squares of coordinates and fixing the value of $r^2$ corresponds to
the definition of
a sphere as the set of points at a fixed distance from the origin.
\\
In our case, by fixing the value of $r^2$ we get the quantum
Euclidean sphere $S^{N-1}_q$
of the corresponding ``radius"\/.  Thus, the quantum Euclidean sphere
$S^{N-1}_q$ is naturally considered as a quantum subspace of the
quantum
Euclidean space $\IR^N_q$ and the algebra $A(S^{N-1}_q)$ of
polynomial functions on it
is a quotient of the  algebra $A(\IR^N_q)$ by the ideal generated by
the relation that fixes the radius. Furthermore, the natural coaction
of $\SO_q(N)$ on
$\IR^N_q$,
\begin{equation}
\delta : A(\IR^N_q) \rightarrow A(\SO_q(N)) \otimes A(\IR^N_q),
\end{equation}
preserves the radius $r^2$, $\delta(r^2) = 1 \otimes r^2$, and yields
a  coaction of
$\SO_q(N)$ on $S^{N-1}_q$.

We can simplify the relations \eqref{old.core} by rescaling one
generator,
$x_0 := (1+q^{- 1})^{-1/2}\xz$. It is also simpler to use rescaled
``partial radii''
which are  related to the $r_i^2$'s by,
\begin{equation*}
r^2_i=(1+q^{-2\rho_i}) s_i ,
\end{equation*}
and are given recursively by
\begin{equation}
\label{pa.ra}
\begin{split}
s_i &:= s_{i-1} + x_i^* x_i = q^{-2} s_{i-1} + x_i x_i^* ,  \\
s_0 &:=x_0^2.
\end{split}
\end{equation}
By using these new elements $s_i$ we can write the
commutation relations of
the generators $\{x_0=x_0^*, x_i, x_i^*, i=1,\dots,n \}$ of the
algebra $A(\IR^N_q)$ as,
\begin{subequations}
\label{core}
\begin{align}
\begin{split}
   x_i x_j &= q x_j x_i ,  \quad 0\leq i<j \leq n, \\
   \qquad x_i^* x_j &= q x_j x_i^*, \quad i\neq j ,
\end{split}
\label{core1} \\
   [x_i, x_i^*] &= (1 - q^{-2}) s_{i-1},
\label{core2}
\end{align}
\end{subequations}
with the understanding that $x_0=0$ if $N=2n$. We see that the
equality of the
two formul\ae\ for $s_i$ in \eqref{pa.ra} is equivalent to the
commutation
relation \eqref{core2}.
\\
Given $N$, the algebras $A(S^{N-1}_q)$ for quantum spheres of
different radii are
isomorphic. We are thus free to normalize the radius to the most
convenient
choice. This is  $s_n=1$ or equivalently $r^2 = 1 + q^{2n-1}$ or $r^2
= 1 +
q^{2n-2}$ according to whether $N$ is odd (even spheres) or $N$ is
even (odd
spheres).
\\
The elements $s_i\in A(S^{2n}_q)$ are self-adjoint and related as
\begin{equation}
\label{order}
0\leq s_0 \leq \dots \leq s_{n-1} \leq s_n = 1.
\end{equation}
     From the commutation relations \eqref{core1} it follows for $i<j$ that
$x_i^*x_ix_j = q^2 x_j x_i^*x_i$; on the other hand $x_j^*x_j x_i =
x_i x_j^*x_j$. With a little induction, we deduce that
\[
s_i x_j =
\begin{cases} q^2 x_j s_i &: i<j \\
\; x_j s_i &: i\geq j,
\end{cases}
\qquad s_i x_j^* =
\begin{cases} q^{-2} x_j^* s_i &: i<j \\
\; x_j^* s_i &: i\geq j
\end{cases}
\]
and that the $s_i$'s are mutually commuting.

Looking at relations \eqref{core}, we see that odd quantum
Euclidean spheres $S^{2n-1}_q$ are the same as the odd quantum
spheres introduced in \cite{VakSoi91} (see also \cite{Wel98}) as
noncommutative
homogeneous spaces of  quantum unitary
groups $\SU_q(n)$.
\\
To our knowledge this simple fact, which was observed during a
conversation
with
F. Bonechi and L. Dabrowski,  has not been presented before. It
extends  the
classical result that odd dimensional spheres are simultaneously
homogeneous spaces  of orthogonal and of unitary groups.

\bigskip
The presentation of the algebra $A(\IR^N_q)$ in \cite{FadResTak89}
involved the square root of the deformation parameter $q$; this must
therefore be positive in that construction.  In our
presentation no square roots are involved and we may take any value
of $q\in\IR$;
however, we shall soon see
that we may restrict the values of $q$ without loss of generality,
due to the
occurrence of
natural isomorphisms.

\bigskip
In $\cite{LanMad01}$ it was shown that the defining relations of the
algebra
$A(S^{N-1}_q)$ of polynomial functions on $S^{N-1}_q$ are equivalent
to
the condition
that a certain matrix over $A(S^{N-1}_q)$ be idempotent. This is also
equivalent to
the condition that another matrix be unipotent, as we shall explain
presently.

First consider the even spheres $S^{2n}_q$ for any integer $n>0$. The
algebra $A(S^{2n}_q)$  is generated by elements $\{ x_0, x_i, x_i^*,
i=1, \dots, n\}$.
\\
We recursively define self-adjoint matrices $u_{(2n)}\in
\Mat_{2^n}(\co\langle x_0, x_i, x_i^*, i=1, \dots, n\rangle)$ for all
$n$ by,
\begin{equation}
\label{uni.even} u_{(2n)} :=
\begin{pmatrix} q^{-1}u_{(2n-2)} & x_n \\ x_n^* & -u_{(2n-2)}
\end{pmatrix}
,\end{equation} with $u_{(0)}=x_0$. The $*$-algebra $A(S^{2n}_q)$
is then defined by the relations that $u_{(2n)}$ is unipotent,
$u_{(2n)}^2=1$, and self-adjoint, $u_{(2n)}^*=u_{(2n)}$. That is, it
is
the quotient of the free $*$-algebra on $2n+1$ generators by these
relations.

The self-adjointness relations merely give that $x_i^*$ is the adjoint
of $x_i$ and $x_0$ is self-adjoint. Unipotency gives a
matrix of $2^{2n}$ relations, although many of these are vacuous or
redundant. These can be deduced inductively from \eqref{uni.even}
which gives,
\begin{equation}
\label{u2recursion}
   u_{(2n)}^2 =
\begin{pmatrix} q^{-2}u_{(2n-2)}^2 + x_nx_n^* & q^{-1}u_{(2n-2)} x_n -
x_n u_{(2n-2)} \\ q^{-1}x_n^* u_{(2n-2)} - u_{(2n-2)} x_n^* &
u_{(2n-2)}^2 + x_n^*x_n
\end{pmatrix}.
\end{equation}
The condition that $u_{(2n)}^2=1$ means in particular that
$u_{(2n)}^2$ is
diagonal with all the diagonal entries equal. Looking at
\eqref{u2recursion}, we
see that the same must be true of $u_{2n-2}^2 \in
\Mat_{2^{n-1}}(A(S^{2n}_q))$,
and so on. Thus, the diagonal relations require that all the diagonal
entries of
(each) $u_{(2j)}^2$ are equal. If this is true for $u_{(2j-2)}^2$,
then the relation
for $u_{(2j)}^2$ is that the same element (the diagonal entry) can be
written in
two different ways. This element is simply $s_j$ and the two ways of
writing it are
those  given in \eqref{pa.ra}. \\
Finally, $u_{(2n)}^2=1$ gives the relation $s_n=1$. \\
The off-diagonal relations are $q^{-1}u_{(2j-2)}x_j = x_ju_{(2j-2)}$
and
$q^{-1}x_j^*u_{(2j-2)} = u_{(2j-2)}x_j^*$ for every $j = 1, \dots ,
n$. Because
the matrix $u_{(2j-2)}$ is constructed linearly from all of the
generators $x_i$
and $x_i^*$ for $i<j$, this conditions are equivalent to the
commutation relations
\eqref{core1}. \\
In summary, up we see that the
defining relations are all obtained from the unipotency condition,
$u_{(2n)}^2=1$.

This presentation of the relations by the unipotency of $u_{(2n)}$ is
the
easiest way to see that there is an
isomorphism $A(S^{2n}_{1/q}) \cong A(S^{2n}_q)$. The substitutions $q
\leftrightarrow q^{-1}$, $x_0 \rightarrow (-q)^{n} x_0$, and $x_i
\rightarrow (-q)^{n-i} x_i^*$ are equivalent to conjugating $u_{(2n)}$
by the antidiagonal matrix
\[
\begin{pmatrix}
& & & & 1\\
& & & 1\\
& & \addots\\
& 1\\
1
\end{pmatrix}
\]
and the result is unipotent and self-adjoint if and only if
$u_{(2n)}$ is; thus there is an isomorphism,
$A(S^{2n}_{1/q})\cong A(S^{2n}_q)$. Because of this, we can assume
that $\abs{q}>1$ without loss of generality.

\bigskip
Now consider the odd spheres $S^{2n-1}_q$ for any integer $n>0$. We
can construct a unipotent $u_{(2n-1)} \in
\Mat_{2^{n}}[A(S^{2n-1}_q)]$, simply by setting
$x_0=0$ in $u_{(2n)}$. Once again, the unipotency condition,
$u_{(2n-1)}^2=1$, is equivalent to the relations defining the algebra
$A(S^{2n-1}_q)$ of polynomial functions on $S^{2n-1}_q$. Again, one
defines self-adjoint elements $s_i\in A(S^{2n-1}_q)$ such  that $s_i =
s_{i-1} + x_i^*x_i = q^{-2} s_{i-1} + x_i x_i^*$ with now $s_0 =
x_0^2 =0$. The
commutation relations are again given by \eqref{core} but now
\eqref{core2} gives in particular that the generator
$x_1$ is normal,
\begin{equation}
x_1x_1^*=x_1^*x_1 \qquad \mathrm{in} \qquad A(S^{2n-1}_q) .
\end{equation}
The previous argument also shows that
$A(S^{2n-1}_q)$ is
the quotient of $A(S^{2n}_q)$ by the ideal generated by
$x_0$; geometrically, this means that $S^{2n-1}_q$ is a
noncommutative subspace of
$S^{2n}_q$.
\\
Because of the  isomorphism $A(S^{2n}_{1/q}) \cong A(S^{2n}_q)$, we
have another isomorphism $A(S^{2n-1}_{1/q}) \cong A(S^{2n-1}_q)$,
  and again we
can  assume that
$\abs{q}>1$ without any loss of generality.

\subsection{Interrelations}
The even sphere algebras each have an involutive
automorphism
\begin{equation}\label{sigma}
\begin{split}
\sigma :\, & A(S^{2n}_q) \to A(S^{2n}_q) \\
& x_0 \mapsto -x_0 , \\
& x_j \mapsto x_j,
\quad  j\neq 0 .
\end{split}
\end{equation}

Obviously, this corresponds to flipping (reflecting) the classical
$S^{2n}$ across the hyperplane $x_0=0$. The coinvariant algebra of
$\sigma$ is the quotient of $A(S^{2n}_q)$ by the ideal generated by
$x_0$, which, as we have noted, is simply
$A(S^{2n-1}_q)$. Geometrically this means that
$S^{2n-1}_q$ is the ``equator'' of $S^{2n}_q$ --- the subspace fixed
by the
flip.

\bigskip
As for odd spheres, they have an action $\rho:\T\to\Aut[A(S^{2n})]$
of the torus group $\T$, defined by
multiplying $x_1$ by a phase
and leaving the  other generators alone,
\begin{equation}\label{rho}
\begin{split}
\rho(\lambda) :\, & A(S^{2n+1}_q) \to A(S^{2n+1}_q) \\
& x_1 \mapsto \lambda x_1 , \\
& x_j \mapsto x_j,
\quad  j\neq 1 .
\end{split}
\end{equation}
The coinvariant algebra is
given by setting $x_1=0$. Now, let  $u'_{(2n+1)}$ be the matrix
obtained by
setting $x_1=0$ and relabeling $x_2$ as $x_1$,
\emph{et cetera}, in the matrix  $u_{(2n+1)}$. Then, $u'_{(2n+1)}$ is
equivalent to
tensoring $u_{(2n-1)}$ with  $\left(\begin{smallmatrix}1 &
0\\ 0&1\end{smallmatrix}\right)$,
\[
u'_{(2n+1)} = u_{(2n-1)} \otimes
\left(\begin{smallmatrix}1 &
0\\ 0&1\end{smallmatrix}\right)
\]
and the result is unipotent if and only if
$u_{(2n-1)}$ is; that is the unipotency of $u'_{(2n+1)}$ yields all
and only
the same relations coming from the unipotency of $u_{(2n-1)}$. This
shows
that
$A(S^{2n-1}_q)$ is the quotient of
$A(S^{2n+1}_q)$ by
the $*$-ideal generated by $x_1$ and $S^{2n-1}_q$ is
the noncommutative
subspace of $S^{2n+1}_q$ fixed by the $\T$-action in \eqref{rho}.

\bigskip
There is also a way of realizing even spheres as noncommutative
subspaces of odd ones. Consider $S^{2n+1}_q$, set $x_1=x_1^*=x_0$ and
relabel
$x_2$ as $x_1$, \emph{et cetera}; let $u''_{(2n+1)}$ be the matrix
obtained from
$u_{(2n+1)}$ with these substitutions. The matrix
$u''_{(2n+1)}$ is the  same as $u_{(2n)}$ in which we substitute
\[
x_0 \rightarrow {\textstyle \begin{pmatrix}0 &
x_0\\ x_0& 0 \end{pmatrix}},
\qquad x_j \rightarrow
{\textstyle\begin{pmatrix} x_{j+1}
& 0 \\ 0 & x_{j+1} \end{pmatrix}}, \quad j\not= 0 \, .
\]
Then the unipotency of $u''_{(2n+1)}$ yields precisely
the same relations coming from the unipotency of $u_{(2n)}$; this
shows
that $A(S^{2n}_q)$ is
the quotient of $A(S^{2n+1}_q)$ by the $*$-ideal generated by
$x_1-x_1^*$.

Thus, every sphere contains a smaller sphere of dimension one less;
by following this tower of inclusions to its  base, we see that every
sphere contains a classical $S^1$, because the circle does not deform.

\section{Structure and Representations}
For
each dimension $N$, we have a one
parameter family of
algebras $A(S^{N-1}_q)$ which, at $q=1$, gives
$A(S^{N-1}_1)=A(S^{N-1})$, the algebra of
polynomial functions on a classical sphere $S^{N-1}$. It is possible
to identify this one-parameter family of algebras to a fixed vector
space and view the product as varying with the parameter. We can then
construct a Poisson bracket on $A(S^{N-1})$ from the first derivative
of the product at the ``classical'' parameter value, $q=1$. The
standard properties of a Poisson bracket (Leibniz and Jacobi
identities) are simple consequences of associativity. The Poisson
bracket is given geometrically by a Poisson bivector $\pi$, an
antisymmetric contravariant $2$-tensor. This in turn determines a
symplectic foliation by the directions in which $\pi$ is
nondegenerate.

In general, given such a one-parameter deformation from a commutative
manifold $\M$ into noncommutative algebras, we can construct a
Poisson bracket on functions.
This Poisson algebra, $A(\M)$ with the commutative product and the
Poisson bracket, describes the deformation to first order.
A deformation is essentially a path through an enormous space of
possible algebras, and the Poisson algebra is just a tangent.
Nevertheless, if the deformation is well behaved the Poisson algebra
does indicate where it is heading. Here are some things
that one can expect.

If $\pi$ vanishes along some subspace $X\subset\M$, then $\pi$
induces a trivial Poisson structure on $X$; i.~e., $X$ is
undeformed to first order. This suggests that it may be undeformed
altogether. If  so, then $X$ will be a classical subspace; that is,
there will be a surjective  homomorphism of the deformed algebra to
the (undeformed) algebra of functions on $X$.

More generally, if $X\subset\M$ is a submanifold such that the
restriction of
$\pi$ to $X$ is tangent to $X$, then  the restriction of functions to
$X$ is a Poisson homomorphism, $A(\M)\to A(X)$; i.~e.,
the deformation respects $X$ to
first order. In this case it may be that some deformation of $X$ is a
``noncommutative subspace'' of the deformation of $\M$. Algebraically
speaking, this  means that the algebra corresponding to $X$ is a
quotient of that corresponding to $\M$. Equivalently, the subalgebra
of functions on
$\M$ vanishing along $X$ correspond to an ideal in the deformation.

Suppose that the symplectic leaves of $\M$ are compact and $Y$ is the
leaf space. Functions which are constant along the symplectic leaves
can be identified with functions on $Y$. In this way, $A(Y)$ is the
center of the Poisson algebra. That is, if $f\in A(Y)$ and $g\in
A(\M)$ then $\{f,g\}=0$. This suggests that the subalgebra
$A(Y)\subset A(\M)$ will be undeformed and will be the center of the
deformed algebra.

More generally, if the symplectic foliation has a Hausdorff leaf
space, $Y$, then $A(Y)$ acts by central multipliers on the Poisson
algebra.
That is for $f$ and $g$ functions on $\M$ and $h$ a
function on $Y$, one has $\{hf,g\} = h\{f,g\}$. This suggests that
$A(Y)$ my undeformed and will act by central multipliers on the
deformation of $A(\M)$.
  If so, then the deformed algebra will
be the  algebra of sections of a bundle of algebras over $Y$. The
fibers will be deformations of the symplectic leaves.

\bigskip
With these ideas in mind, consider some of the properties of
the deformed spheres $S^{N-1}_q$.

We have seen that the $S^{2n-1}_q$ noncommutative subspace of
$S^{2n}_q$
corresponds to the equator, the $S^{2n-1}\subset S^{2n}$ where
$x_0=0$ and the
Poisson bivector  on $S^{2n}$ is degenerate. On the remaining
$S^{2n}\smallsetminus S^{2n-1}$, the Poisson bivector is
nondegenerate.
So, topologically, we have a union of two copies of symplectic
$\R^{2n}$.
\\
The manifold $\R^{2n}$ has a unique symplectic structure, modulo
isomorphism. This symplectic structure corresponds to an essentially
unique deformation. If we complete to a \cs-algebra, then the
deformation
of $\C_0(\R^{2n})$ will be the algebra, $\K$, of compact operators on
a
countably infinite-dimensional Hilbert space.
\\
Then, the kernel of the quotient map $A(S^{2n}_q)\to A(S^{2n-1}_q)$
should be a deformation of the subalgebra of functions on $S^{2n}$
which
vanish at the equator. If we complete to \cs-algebras, this should
give us the  direct sum of two copies of $\K$, one for each
hemisphere.
Thus we expect that the
\cs-algebra $\C(S^{2n}_q)$ will be an extension:
\begin{equation}\label{ext.even}
0 \to \K\oplus\K \to \C(S^{2n}_q) \to \C(S^{2n-1}_q)
\to 0
.\end{equation}

In odd dimensions, the Poisson bivector is necessarily degenerate.
However, the
$S^{2n-1}_q$ noncommutative subspace of $S^{2n+1}_q$ corresponds
classically to the Poisson bivector being more degenerate on
$S^{2n-1}\subset  S^{2n+1}$. It is of rank $2n$ at most points, but of
rank $2n-2$ (or less) along
$S^{2n-1}$. The complement
$S^{2n+1}\smallsetminus S^{2n-1}$ has a symplectic foliation by
$2n$ dimensional leaves which is invariant under the $\T$ action;
the simplest possibility is that this corresponds to the product in
the
identification
\[
S^{2n+1}\smallsetminus S^{2n-1} \cong S^1\times \R^{2n}
.\]
If we complete to \cs-algebras, then the deformation of this should be
$\C(S^1)\otimes\K$. The kernel of the quotient map $A(S^{2n+1}_q)\to
A(S^{2n-1}_q)$ should be this deformation, so we expect another
extension,
\begin{equation}
\label{ext.odd} 0 \to \C(S^1)\otimes\K \to \C(S^{2n+1}_q)\to
\C(S^{2n-1}_q) \to 0
.\end{equation}

These expectations are true. As we
have mentioned, the odd dimensional spheres we are considering are
equivalent to
the ``unitary" odd quantum spheres of Vaksman and Soibelman
\cite{VakSoi91}. In
\cite{HonSzy01} Hong and \Szymanski\ obtained the \cs-algebras
$\C(S^{2n+1}_q)$ as Cuntz-Krieger algebras of suitable graphs. From
this
construction they derived the extension \eqref{ext.odd}. They
also considered even spheres, defined as quotients of odd ones  by
the ideal
generated by $x_1 - x_1^*$. These are thus isomorphic to the even
spheres we are
considering here. They also obtained these as Cuntz-Krieger algebras
and derived
the extension \eqref{ext.even}. However, as explicitly stated in the
introduction to \cite{HonSzy01}, they were unable to realize even
spheres as
quantum homogeneous spaces of quantum orthogonal groups, thus also
failing to
realize that ``unitary" and ``orthogonal" odd quantum spheres are the
same.

\bigskip
Representations of the odd dimensional spheres were constructed in
\cite{VakSoi91}. The primitive spectra of all these
spheres were
compute in \cite{HonSzy01},  which amounts to a classification of
representations. The
representations for
quantum Euclidean spheres have also been constructed in \cite{Fio95}
by thinking
of them as quotient algebras of quantum Euclidean planes.
\\
By using the properties we have just discussed, we shall present
a clearer derivation of the representations. Indeed, the structure of
the representations can be anticipated from the construction of
$S^{N-1}_q$
via the extensions \eqref{ext.even} and \eqref{ext.odd} and by
remembering that an
irreducible representation $\psi$ can be partially characterized by its
kernel. Moreover, an irreducible representation of a \cs-algebra restricts either to an irreducible or a trivial representation of any ideal; and conversely, an irreducible representation of an ideal extends to an irreducible representation of the \cs-algebra (see for instance \cite{Dav96}).

For an even sphere $S^{2n}_q$, the kernel of an irreducible
representation $\psi$ will
contain one or both of the copies of $\K\subset \C(S^{2n}_q)$.
If $\K\oplus\K\subseteq\ker\psi$, then $\psi$ factors through
$\C(S^{2n-1}_q)$ and is given by a representation of that algebra.
If one copy of $\K$ is not in $\ker\psi$, then $\psi$
restricts to a representation of this $\K$. However, $\K$ has only one
irreducible representation. Since $\K$ is an ideal in $\C(S^{2n}_q)$,
the unique irreducible representation of
$\K$ uniquely extends to a representation of $\C(S^{2n}_q)$ (with the
other copy of $\K$
in its kernel).
\\
Thus, we expect the irreducible representations of $S^{2n}_q$ (up to
isomorphism) to be:

1. all irreducible representations of $S^{2n-1}_q$

2. a unique representation with kernel the second copy of $\K$

3. a unique representation with kernel the first copy of $\K$.

  From the extension \eqref{ext.even} we expect that the generator
$x_0$ is a
self-adjoint  element of $\K\oplus\K\subset \C(S^{2n}_q)$ and it
should have almost
discrete, real spectrum: it will therefore be a convenient tool for
decomposing
the Hilbert space in  a representation.

\bigskip
Similarly, from the construction of $S^{2n+1}_q$ by the extension
\eqref{ext.odd}, one can anticipate the structure of its
representations.
Firstly, if $\C(S^1)\otimes\K\subseteq\ker\psi$, then $\psi$ factors
through
$\C(S^{2n-1}_q)$ and is really a representation of $S^{2n-1}_q$.

Otherwise, $\psi$ restricts to an irreducible representation of
$\C(S^1)\otimes\K$. This factorizes as the tensor product of an
irreducible representation of $\C(S^1)$ with one of $\K$. The
irreducible representations of $\C(S^1)$ are simply given by the
points of
$S^1$, and as we have
mentioned, $\K$ has a unique irreducible representation. The
representations of
$\C(S^1)\otimes\K$ are thus classified by the points of $S^1$. These
representations extend  uniquely from the ideal $\C(S^1)\otimes\K$ to
the whole algebra $\C(S^{2n+1}_q)$.
\\
Thus, we expect the irreducible representations of $S^{2n+1}_q$ (up to
isomorphism) to be:

1. all irreducible representations of $S^{2n-1}_q$

2. a family of representations parameterized by $S^1$.

In our construction of the representations, a simple identity regarding the spectra of operators will be especially useful (see, for instance \cite{Sak98}). If $x$ is an element of any
\cs-algebra, then
\begin{equation}
\label{spec.id}
\{0\}\cup \Spec x^*x = \{0\}\cup \Spec xx^*
.\end{equation}

%\subsection{A Spectral Lemma}

%This is easy to see if $x$ is a compact operator on a Hilbert space.
%In that
%case these spectra consist of eigenvalues (except possibly $0$).
%Suppose that
%$\lambda\in\co$ is a nonzero eigenvalue of $x^*x$ and $\ket\varphi$ a
%corresponding eigenvector. Then $x\ket\varphi \neq0$ and
%\[
%(xx^*)x\ket\varphi = x(x^*x)\ket\varphi = \lambda x\ket\varphi
%\]
%so $x\ket\varphi$ is an eigenvector of $xx^*$ with eigenvalue
%$\lambda$.
%Therefore
%\[
%\Spec x^*x \subseteq \{0\}\cup \Spec xx^*
%\]
%   and with the symmetric result, \eqref{spec.id} follows.

%We cannot assume that $x$ is compact in general. Indeed,
%\eqref{spec.id} is
%useful in proving that an operator is compact. In general, suppose
%that
%$\lambda\neq0\in\co$ and $\lambda\not\in\Spec xx^*$. This implies
%that the
%inverse, $(xx^*-\lambda)^{-1}$ exists. Let $a=x^*(xx^*-\lambda)^{-1}x
%- 1$. By
%direct computation,
%\[
%\begin{split}
%(x^*x-\lambda)a
%&= (x^*x-\lambda)[x^*(xx^*-\lambda)^{-1}x - 1] \\
%&= x^*(xx^*-\lambda)(xx^*-\lambda)^{-1}x - (x^*x-\lambda) \\
%&= \lambda
%\end{split}
%\]
%and so $\lambda^{-1}a$ is a left (and similarly, right) inverse of
%$x^*x-
%\lambda$. Therefore $\lambda\not\in\Spec x^*x$. Therefore $\Spec x^*x
%\subseteq
%\{0\}\cup \Spec xx^*$ and \eqref{spec.id} follows in general.

\subsection{Even Sphere Representations}
To illustrate the general structure we shall start by describing the  lowest
dimensional case, namely $S^2_q$. This is isomorphic to the so-called equator
sphere of Podle{\`s} \cite{Pod87}. For this  sphere, the 
representations were also
constructed in \cite{MasNakWat91} in a way  close to the one presented here.

Let us then consider the sphere $S^2_q$.
\\
As we have discussed, we expect that $x_0$ is a compact operator (in
some
faithful representation) and thus has an
almost discrete, real spectrum. The relation $x_0 x_1=qx_1x_0$
suggests
that $x_1$ and
$x_1^*$ shift the eigenvalues of $x_0$. However, we cannot assume
\emph{a priori} that $x_0$ has eigenvalues, let alone that
eigenvectors
form a complete basis of the Hilbert space.
The relation $1 = x_0^2 + x_1^*x_1 = q^{-2} x_0^2 + x_1 x_1^* $ shows
that
$x_0^2 \leq 1$ and  thus
$\norm{x_0}
\leq 1$. As $x_0$ is self-adjoint, this shows that $\Spec x_0
\subseteq [-1,1]$. By \eqref{spec.id} we have also,
\begin{align*}
\{0\}\cup \Spec x_1^*x_1 &= \{0\}\cup \Spec x_1x_1^* \\
\{0\}\cup \Spec (1-x_0^2) &= \{0\}\cup \Spec (1- q^{-2}x_0^2) \\
\{1\}\cup \Spec x_0^2 &= \{1\} \cup q^{-2} \Spec x_0^2
.\end{align*}
Because we have assumed that $\abs{q}\geq1$, the only subsets of
$[0,1]$ that
satisfy this condition are $\{0\}$ and
$\{0,q^{-2k}\mid k=0,1,\dots\}$.
\\
If $x_0\neq 0 \in \C(S^2_q)$ then $\Spec x_0^2$ is the latter set. We
cannot
simply assume that $x_0\neq 0$, since not every $*$-algebra is a
subalgebra of a
\cs-algebra; however, our explicit representations will show that
that is the
case here.

Now let $\Hi$ be a separable Hilbert space and suppose that $\psi :
A(S^2_q) \to
\Li(\Hi)$ is an irreducible $*$-representation.

If $\psi(x_0)= 0$ then $1 = \psi(x_1)\psi(x_1)^*
=\psi(x_1)^*\psi(x_1)$.
Thus
$\psi(x_1)$ is unitary, and by the assumption of
irreducibility, it is a number $\lambda\in\co$, $\abs\lambda = 1$. So,
$\Hi=\co$ and the representation is $\psi^{(1)}_\lambda$ defined by,
\begin{equation}\label{reps1}
\begin{split}
\psi^{(1)}_\lambda(x_0) &=0,  \\
\psi^{(1)}_\lambda(x_1) &=\lambda , \qquad \lambda\in S^1 .
\end{split}
\end{equation}
Thus we have an $S^1$ worth of representations with $x_0$ in the
kernel.

If $\psi(x_0)\neq 0$, then $1 \in \Spec x_0^2$; it is an isolated point in the spectrum and therefore an eigenvalue. For some sign
$\pm$ there exists a unit vector $\ket0 \in\Hi$ such that
$\psi(x_0)\ket0 =
\pm \ket 0$. So, for $k=0,1,\dots$, $\psi(x_1^*)^k\ket 0$ is an
eigenvector as well, because
\[
\psi(x_0)  \psi(x_1^*)^k\ket 0  = q^{-k} \psi({x_1^*}^k x_0) \ket 0  =
\pm q^{-k}
\psi(x_1^*)^k\ket 0
.\] By normalizing, we obtain a sequence of unit
eigenvectors, recursively defined by
\[
\ket k := (1-q^{-2k})^{-1/2} \psi(x_1^*)\ket{k-1}
.\]
We have thus two representations $\psi^{(2)}_+$ and $\psi^{(2)}_-$,
and direct computation shows that
\begin{equation}\label{reps2}
\begin{split}
       \psi^{(2)}_\pm(x_0) \ket k &= \pm
       q^{-k} \ket k,  \\
\psi^{(2)}_\pm(x_1) \ket k &= (1 - q^{-2k})^{1/2}
\ket{k-1} ,  \\
\psi^{(2)}_\pm(x_1^*) \ket k &= (1 - q^{-2(k+1)})^{1/2} \ket{k+1} .
\end{split}
\end{equation}
The eigenvectors $\{\ket k\mid k=0,1,\dots\}$
are mutually orthogonal because they have distinct eigenvalues, and by
the assumption of irreducibility form a basis for $\Hi$.
\\
Notice that any power of $\psi^{(2)}_\pm(x_0)$ is a trace class
operator, while this is not the case for the operators
$\psi^{(2)}_{\pm}(x_1)$ and $\psi^{(2)}_{\pm}(x_1^*)$ nor for any of
their powers.
\\
Note also that the representations \eqref{reps2} are related by the
automorphism
$\sigma$ in
\eqref{sigma}, as
\begin{equation}
\psi^{(2)}_\pm\circ\sigma = \psi^{(2)}_\mp.
\end{equation}

If we set a value of $q$ with $\abs{q}<1$ in \eqref{reps2}, the
operators would be
unbounded. This is the reason for assuming that $\abs{q}>1$. The
assumption
was used in computing
$\Spec x_0$. Not only is $\norm{x_0}\leq 1$, but by a similar
calculation
$\norm{x_0}\leq \abs{q}$. Which bound is more relevant obviously
depends on
whether $q$ is greater or less than $1$.
\\
For $\abs{q}<1$ the appropriate formul\ae\ for the representations
can be obtained from \eqref{reps2}
by replacing the index $k$ with $-k-1$. As a consequence, the role of
$x_1$
and $x_1^*$ as lowering and raising operators is exchanged.

\bigskip
Let us then turn to the general even spheres $S^{2n}_q$. \\
The structure of the representations is similar to that for $S^2_q$,
but more
complicated. The element $x_0$ is no longer
sufficient to completely decompose the Hilbert space of the
representation and we need to
use all the commuting self-adjoint elements $s_i\in A(S^{2n}_q)$
defined in
\eqref{pa.ra}.

Using the formul\ae\ for $1 = s_n = s_{n-1} + x_n^*x_n
= q^{-2}s_{n-1} + x_n x_n^*$
   and \eqref{spec.id} we get,
\begin{align*}
\{0\}\cup \Spec x_n^*x_n &= \{0\} \cup \Spec x_n x_n^* \\
\{0\}\cup \Spec (1-s_{n-1}) &= \{0\} \cup \Spec (1-q^{-2}s_{n-1}) \\
\{1\} \cup \Spec s_{n-1} &= \{1\} \cup q^{-2} \Spec s_{n-1}
.\end{align*}
Therefore either $s_{n-1} =0\in\C(S^{2n}_q)$ or $\Spec s_{n-1} =
\{0,q^{-2k}\mid k=0,1,\dots\}$.

Now suppose that $\psi : A(S^{2n}_q) \to \Li(\Hi)$ is an irreducible
$*$-representation.

If $\psi(x_0) =0$, then $\psi$ factors through $A(S^{2n-1}_q)$. Thus
$\psi$ is an irreducible representation of $A(S^{2n-1}_q)$; these
will be
discussed later.

If $\psi(x_0) \neq 0$, then $\psi(s_0)\neq 0$, and by the relations
\eqref{order}, all of the $\psi(s_i)$'s are nonzero. So, $\Spec
\psi(s_{n-1}) =
\{0,q^{-2k}\mid k=0,1,\dots\}$ and in particular, $1$ is an
eigenvalue of
$\psi(s_{n-1})$. Because $s_{n-1}$ commutes with all the generators
except $x_n$
and $x_n^*$, these generators preserve the corresponding eigenspace.
The same argument as for $s_{n-1}$ shows that the restriction of
$\psi(s_{n-2})$ to this
eigenspace has  the same spectrum; in particular, $1$ is an
eigenvalue. There is a
simultaneous  eigenspace of $\psi(s_{n-1})$ and $\psi(s_{n-2})$ with
eigenvalue $1$ for
both. Proceeding in this way, we find that there is a simultaneous
eigenspace with
eigenvalue $1$ for all the $\psi(s_i)$'s. That is, there must exist a
unit vector
$\ket{0,\dots,0}\in \Hi$ such that $\psi(s_i) \ket{0,\dots,0} =
\ket{0,\dots,0}$ for all
$i$ and $\psi(x_0) \ket{0,\dots,0} = \pm \ket{0,\dots,0}$. We
construct more unit
vectors by letting $\ket{k_0,\dots,k_{n-1}}$ be
\[
\psi(x_1^*)^{k_0}\dots \psi(x_n^*)^{k_{n-1}} \ket{0,\dots,0}
\]
modulo a positive normalizing factor. Using the commutation
relations between
$s_i$'s and $x_j^*$'s we get that,
\[
\psi(s_i) \ket{k_0,\dots,k_{n-1}} = q^{-2(k_i+\dots+k_{n-1})}
\ket{k_0,\dots,k_{n-1}}
.\]
In summary, $x_i$ lowers $k_{i-1}$, $x_i^*$ raises $k_{i-1}$, and
$s_i$ measures the sum $k_i+\dots+k_{n-1}$.
\\
The correct normalizing factors can be determined from,
\begin{align*}
\braket{k_0,\dots,k_{n-1}}{\psi(x_ix_i^*)}{k_0,\dots,k_{n-1}} &=
\braket{k_0,\dots,k_{n-1}}{\psi(s_i-s_{i-1})}{k_0,\dots,k_{n-1}} \\ &=
(1-q^{-2k_{i-1}})q^{-2(k_i+\dots+k_{n-1})}
.\end{align*}
Thus we get two representations $\psi^{(2n)}_\pm$ defined by,
\begin{align}
\label{rep.even}
\psi^{(2n)}_\pm(x_0) \ket{k_0,\dots,k_{n-1}} &= \pm
q^{-(k_0+\dots+k_{n-1})}
\ket{k_0,\dots,k_{n-1}}  \\
\psi^{(2n)}_\pm(x_i) \ket{k_0,\dots,k_{n-1}} &= (1-
q^{-2k_{i-1}})^{1/2}
q^{-(k_i+\dots+k_{n-1})} \ket{k_0,\dots,k_{i-1}-1,\dots,k_{n-1}} \nn
\\
\psi^{(2n)}_\pm(x_i^*) \ket{k_0,\dots,k_{n-1}} &= (1-
q^{-2(k_{i-1}+1)})^{1/2}
q^{-(k_i+\dots+k_{n-1})} \ket{k_0,\dots,k_{i-1}+1,\dots,k_{n-1}} \nn
\end{align}
with $i=1,\dots,n$. The assumption of irreducibility implies that the
collection of
vectors
$\{ \ket{k_0,\dots,k_{n-1}} , k_i \geq 0 \}$ constitute a complete
basis for $\Hi$.
\\
As before, the two representations \eqref{rep.even} are related by
the automorphism $\sigma$, as
\begin{equation}
\psi^{(2n)}_\pm\circ\sigma = \psi^{(2n)}_\mp.
\end{equation}
Again the formul\ae\ \eqref{rep.even} for the representations are for
$\abs{q}>1$; and again the representations for $\abs{q}<1$ can be
obtained
by replacing all indices $k_i$ with $-k_i-1$.

In all of the irreducible representations of $A(S^{2n}_q)$, the
representative of $x_0$ is compact; in fact it is trace class. We can
deduce from
this that the \cs-ideal generated by $\psi^{(2n)}_\pm(x_0)$ in
$\C(S^{2n}_q)$ is
isomorphic to  $\K(\Hi)$, the ideal of all compact operators on
$\Hi$. By using the
continuous functional calculus, we can apply any function
$f\in\C[-1,1]$ to $x_0$. If $f$ is supported on $[0,1]$, then $f(x_0)
\in \ker \psi^{(2n)}_-$. Likewise if $f$ is supported in $[-1,0]$,
then
$f(x_0) \in \ker \psi^{(2n)}_+$.
     From this we deduce that the \cs-ideal generated by $x_0$ in
$\C(S^{2n}_q)$ is $\K\oplus\K$. One copy of $\K$ is $\ker
\psi^{(2n)}_+$; the other is $\ker \psi^{(2n)}_-$. This shows that
the
extension
\eqref{ext.even} is correct.

\subsection{Odd Sphere Representations}
Again, to illustrate the general strategy we shall work out in detail 
the simplest
case, that of the sphere $S^3_q$. This can be identified with the underlying
noncommutative space of the quantum group $\SU_q(2)$  and as such the
representations of the algebra are well known \cite{Wor87}.

The generators $\{x_i, x_i^* \mid
i=1,2 \}$ of the algebra $A(S^3_q)$ satisfy the commutation relations
$x_1 x_2 = q x_2 x_1$, $x_i^* x_j = q x_j x_i^*$,$i\neq j$,
$[x_1, x_1^*] = 0$, and $[x_2, x_2^*] = (1-q^{-2})x_1 x_1^*$.
Furthermore,
there is the sphere relation $1 = x_2^* x_2 + x_1^* x_1 =  x_2 x_2^* +
q^{-2} x_1 x_1^*$.
\\
The normal generator $x_1$ plays much the same role for the
representations of
$S^3_q$ that $x_0$ does for those of $S^2_q$. The sphere relation
shows that
$\norm{x_1}\leq 1$ and
\begin{align*}
\{0\} \cup \Spec x_2^*x_2 &= \{0\} \cup \Spec x_2x_2^* \\
\{0\} \cup \Spec (1 - x_1^*x_1) &= \{0\} \cup \Spec (1 - q^{-2} x_1
x_1^*) \\
\{1\} \cup \Spec x_1^*x_1 &= \{1\} \cup q^{-2} \Spec x_1^*x_1
,\end{align*}
which shows that either $x_1=0$ or $\Spec x_1^*x_1 =
\{0,q^{-2k}\mid k=0,1,\dots\}$.

Let $\psi : A(S^3_q) \to
\Li(\Hi)$ be an irreducible $*$-representation.

If $\psi(x_1)=0$ then
the relations reduce to $1 =
\psi(x_2)\psi(x_2)^* =
\psi(x_2)^*\psi(x_2)$. Thus $\psi(x_2)$ is
unitary and by the
assumption of irreducibility, it is a scalar,
$\psi(x_2)=\lambda \in
\co$ with
$\abs\lambda=1$. Thus, as before, we have an
$S^1$ of representations of this kind.

If $\psi(x_1)\neq0$, then
$1\in \Spec \psi(x_1^*x_1)$ and is an isolated point in the spectrum. Thus, there exists
a unit vector $\ket0
\in \Hi$ such that $\psi(x_1^*x_1) \ket 0 = \ket
0$, and by the
assumption of irreducibility, there is some
$\lambda\in\co$
with
$\abs\lambda=1$ such that $\psi(x_1)\ket0 =
\lambda\ket0$. We see then that
$\psi(x_2^*)^k\ket0$ is an
eigenvector
\[
\psi(x_1) \psi(x_2^*)^k\ket0  = q^{-k}
\psi({x_2^*}^kx_1)\ket0  =
\lambda q^{-k}
\psi(x_2^*)^k\ket0
.\]
By normalizing, we get a sequence of unit
eigenvectors
recursively defined by
\[
\ket k := (1-q^{-2k})^{-1/2}
\psi(x_2^*) \ket{k-1}
.\]
A family of
representations $\psi^{(3)}_\lambda$, $\lambda \in S^1$, is then
defined by
\begin{align}\label{reps3}
\psi^{(3)}_{\lambda}(x_1) \ket{k} &=
\lambda q^{-k} \ket{k}
, \nn \\
\psi^{(3)}_{\lambda}(x_1^*)
\ket{k} &=  \bar{\lambda} q^{-k} \ket{k}
, \nn
\\
\psi^{(3)}_{\lambda}(x_2) \ket{k} &=  (1 - q^{-2k})^{1/2}
\ket{k-1}
, \nn \\
\psi^{(3)}_{\lambda}(x_2^*) \ket{k} &=  (1 -
q^{-2(k+1)})^{1/2}
\ket{k+1} .
\end{align}
We notice that
any power of $\psi^{(3)}_{\lambda}(x_1)$
or
$\psi^{(3)}_{\lambda}(x_1^*)$ is a trace class operator, while
this is
not the case for the operators $\psi^{(3)}_{\lambda}(x_2)$ and
$\psi^{(3)}_{\lambda}(x_2^*)$ nor for any of their powers.

\bigskip

Consider the general odd spheres $S^{2n+1}_q$ and
let $\psi:A(S^{2n+1}_q)\to\Li(\Hi)$ be an irreducible
representation.

If $\psi(x_1) =0$ then $\psi$ factors through $A(S^{2n-1}_q)$ and is
an
irreducible representation of that algebra.

If $\psi(x_1)\neq 0$ then $\psi(s_1)\neq 0, \psi(s_2)\neq 0, $
\emph{et
cetera}. By the same
arguments as  for $S^{2n}_q$, there must exist a simultaneous
eigenspace with eigenvalue
$1$  for all of $s_1,\dots s_n$. By the assumption of irreducibility,
this
eigenspace is
$1$-dimensional.
Let
$\ket{0,\dots,0} \in \Hi$ be a unit vector in this eigenspace. Then
$s_i
\ket{0,\dots,0}
=
\ket{0,\dots,0}$ for $i=1,\dots,n$. The restriction of $\psi(x_1)$ to
this
subspace is unitary and thus for some $\lambda\in\co$ with
$\abs\lambda=1$,
$\psi(x_1) \ket{0,\dots,0} = \lambda \ket{0,\dots,0}$.
   We can construct more simultaneous eigenvectors of the $s_i$'s.
Define $\ket{k_1,\dots,k_n}$ to be
\[
\psi(x_2)^{k_1}\dots \psi(x_{n+1})^{k_n}
\ket{0,\dots,0}
\]
modulo a positive normalizing constant. Then
\[
\psi(x_1) \ket{k_1,\dots,k_n} = \lambda \ket{k_1,\dots,k_n}
\]
   and
\[
\psi(s_i) \ket{k_1,\dots,k_n} = q^{-2(k_i+\dots + k_n)}
\ket{k_1,\dots,k_n}.
\]
Working out the normalization, this gives a family of representations
$\psi^{(2n+1)}_\lambda$ by,
\begin{equation}
\begin{split}
\label{rep.odd}
\raisetag{10ex}
\psi^{(2n+1)}_{\lambda}(x_1)
\ket{k_1,\dots,k_n} &= \lambda
q^{-(k_1+\dots+k_n)}
\ket{k_1,\dots,k_n},   \\
\psi^{(2n+1)}_{\lambda}(x_1^*)
\ket{k_1,\dots,k_n} &=  \bar{\lambda}
q^{-(k_1+\dots+k_n)}
\ket{k_1,\dots,k_n},   \\
\psi^{(2n+1)}_{\lambda}(x_i)
\ket{k_1,\dots,k_n} &=  (1-
q^{-2k_{i-1}})^{1/2}
q^{-(k_i+\dots+k_n)}
\ket{k_1,\dots,k_{i-1}-1,\dots,k_n},  \\
\psi^{(2n+1)}_{\lambda}(x_i^*) \ket{k_1,\dots,k_n} &=
(1-q^{-2(k_{i-1}+1)})^{1/2} q^{-(k_i+\dots+k_n)}
\ket{k_1,\dots,k_{i-1}+1,\dots,k_n} ,
\end{split}
\end{equation}
for
$i=2,\dots,n+1$. The assumption of
irreducibility implies that the
vectors $\{ \ket{k_1,\dots,k_n} , k_i \geq 0 \}$ form an orthonormal
basis of $\Hi$.

As for the even case, the formul\ae\ \eqref{rep.odd} give bounded
operators only for $\abs{q}>1$; and as before, the representations for
$\abs{q}<1$ can be obtained by replacing all indices $k_i$ with
$-k_i-1$.

Again, as in the even case, we can verify  that
$\psi^{(2n+1)}_\lambda(x_1)$ is compact (indeed, trace class) and
that the  ideal
generated by $\psi^{(2n+1)}_\lambda(x_1)$ in the \cs-algebra
completion of the  image
$\psi^{(2n+1)}_\lambda(A(S^{2n+1}_q))$ is  $\K(\Hi)$. The
representations
$\psi^{(2n+1)}_\lambda$ can  be assembled into a single
representation by adjointable
operators on  a Hilbert $\C(S^1)$-module. With this we can verify
that the ideal
generated by $x_1$ in $\C(S^{2n+1}_q)$ is $\C(S^1)\otimes\K$ and this
fact verifies the
extension \eqref{ext.odd}.

\bigskip Putting together  the results for even and odd spheres, we
get a complete
picture of the set of  irreducible representations of all these
spheres.

For the odd spheres $S^{2n+1}_q$, the set  of irreducible
representations (or
equivalently, the primitive spectrum of $\C(S^{2n+1}_q)$) is indexed
by
the union of $n+1$ copies of $S^1$.  These run from the
representations
$\psi^{(2n+1)}_\lambda$ of $S^{2n+1}_q$ given in \eqref{rep.odd} down
to the one dimensional
representations $\psi^{(1)}_\lambda$ that factor through  $\C(S^1)$.

For the even spheres $S^{2n}_q$, the set of irreducible
representations (or
equivalently, the primitive spectrum of $\C(S^{2n}_q)$) is indexed by
the union of
$n$ copies of $S^1$ and $2$ points. The isolated points correspond to
the $2$
representations $\psi^{(2n)}_\pm$ specific to $S^{2n}_q$ and given in
\eqref{rep.even};
the circles  correspond to  representations $\psi^{(2m+1)}_\lambda$
coming from lower
odd  dimensional spheres, down to $S^1$.

\section{$K$-homology and $K$-theory}
We are now ready to study the $K$-homology and $K$-theory of the
quantum Euclidean
spheres $S^{N-1}_q$. The $K$-homology classes will be given by
Fredholm modules
using the representations constructed previously while the $K$-theory
classes
will
be given by means of suitable idempotents and unitaries.
\\
In fact, in order to compute the  pairing of $K$-theory with
$K$-homology, it is
more convenient to first compute the Chern characters and then use
the pairing
between
cyclic homology and cohomology \cite{Con94}. Thus, together with the
generators
of
$K$-theory and $K$-homology we shall also construct the associated
Chern
characters in the cyclic homology $HC_*[A(S^{N-1}_q)]$ and cyclic
cohomology $HC^*[A(S^{N-1}_q)]$ respectively.

It is worth recalling the $K$-theory
and homology of the
classical spheres.
\\
For an even dimensional sphere $S^{2n}$, the
groups are
\begin{gather*} K^0(S^{2n}) \cong \Z^2 ,\qquad K^1(S^{2n}) =0 , \\
K_0(S^{2n})
\cong \Z^2 ,\qquad K_1(S^{2n}) =0.
\end{gather*} One generator of the $K$-theory $[1]\in K^0(S^{2n})$ is
given by
the trivial $1$-dimensional bundle. The other generator of
$K^0(S^{2n})$ is the
left handed  spinor bundle. One $K$-homology generator $[\varepsilon]\in
K_0(S^{2n})$ is  ``trivial'' and is the push-forward of the generator
of $K_0(*)\cong\Z$ by the inclusion
$\iota:*\into S^{2n}$ of a point (any point) into the sphere. The
other generator,
$[\mu]\in  K_0(S^{2n})$, is the $K$-orientation of
$S^{2n}$ given by
its
structure as a spin manifold \cite{Con94}.

For an odd dimensional sphere, the groups are
\begin{gather*} K^0(S^{2n+1}) \cong \Z ,\qquad K^1(S^{2n+1}) \cong \Z
,
\\ K_0(S^{2n+1}) \cong \Z ,\qquad K_1(S^{2n+1}) \cong \Z .
\end{gather*} The generator $[1]\in K^0(S^{2n+1})$ is the
trivial
$1$-dimensional bundle. The generator of $K^1(S^{2n+1})$ is
a
nontrivial unitary matrix-valued function on $S^{2n+1}$. The
generator
$[\varepsilon]\in K_0(S^{2n+1})$  is again the ``trivial''
element given
by the inclusion of a point. The generator $[\mu]\in
K_1(S^{2n+1})$ is the $K$-orientation of $S^{2n+1}$ given by its
structure as a spin manifold \cite{Con94}.

There is a natural pairing between $K$-homology and $K$-theory. If we
pair
$[\varepsilon]$ with a vector bundle we get the rank of the vector
bundle, i.~e.\ the dimension of its fibers. If we pair $[\mu]$ with a
vector bundle it gives the ``degree'' of the bundle, a measure of its
nontriviality. Similarly, pairing with $[\mu]$ measures the
nontriviality of a unitary.

\bigskip
The $K$-theory and $K$-homology of the quantum Euclidean spheres are
isomorphic to that of the classical spheres; that is, for any $N$ and
$q$, $K_*[\C(S^{N-1}_q)] \cong K^*(S^{N-1})$ and $K^*[\C(S^{N-1}_q)]
\cong K_*(S^{N-1})$.

In the case of $K$-theory, this was proven by Hong and \Szymanski\ in
\cite{HonSzy01} using their construction of the \cs-algebras as
Cuntz-Krieger algebras of graphs. The groups $K_0$ and $K_1$ were  given as the
cokernel and the kernel respectively, of a matrix canonically 
associated with
the graph.  The result for $K$-homology can be proven using the same techniques
\cite{Cun84,RaeSzy00}: the groups $K^0$ and $K^1$ are now given as 
the kernel and
the cokernel respectively, of the transposed matrix.
The $K$-theory and the $K$-homology for the particular case of $S^2_q$ (in fact
for all Podle{\`s} spheres $S^2_{qc}$) was worked out in \cite{MasNakWat91})
while for $S^3\cong SU_q(2)$ it was spelled out in \cite{MasNakWat90}.

\subsection{$K$-homology}
\label{se:khom}
Because the  $K$-homology of these deformed spheres is isomorphic to
the
$K$-hom\-o\-logy of the ordinary spheres, we need to construct  two
independent
generators. First consider the ``trivial'' generator of
$K^0[\C(S^{N-1}_q)]$.  This can be constructed in a manner closely
analogous to the  undeformed case.

As we have just described, the trivial generator of $K_0(S^{N-1})$
is  the
image of the
generator of the $K$-homology of a point by the  functorial map
$K_*(\iota) : K_0(*) \to K_0(S^{N-1})$, where $\iota : *\into
S^{N-1}$
is the inclusion
of a point into the sphere. The  quantum Euclidean spheres do not
have as many points,
but they do have some.
We have seen that the relations among the various spheres always
include
a  homomorphism
$A(S_q^{N-1}) \to A(S^1)$. Equivalently, every $S^{N-1}_q$ has  a
circle $S^1$ as a
classical subspace;  thus for  every
$\lambda\in S^1$ there is a point, i.~e., the  homomorphism
$\psi^{(1)}_\lambda :
\C(S^{N-1}_q) \to \co$.

We can construct an element $[\varepsilon_\lambda] \in K^0[\C(S_q^{N-1})]$ by pulling back the generator of $K^0(\co)$ by $\psi^{(1)}_\lambda$. This construction factors through $K_0(S^1)$. Because $S^1$ is path connected, the points of $S^1$ all define homotopic (and hence $K$-homologous) Fredholm modules. Thus there is a single $K$-homology class $[\varepsilon_\lambda] \in K^0[\C(S_q^{N-1})]$, independent of $\lambda\in S^1$. 
\\
The canonical generator of $K^0(\co)$ is  given by the following
Fredholm module: The
Hilbert space is $\co$; the grading operator is $\gamma=1$; the
representation is the
obvious representation of $\co$ on $\co$; the Fredholm operator is
$0$. If we  pull this
back to $K^0[\C(S^{N-1}_q)]$ using
$\psi^{(1)}_\lambda$, then  the Fredholm module
$\varepsilon_\lambda$ is given in the
same way  but with
$\psi^{(1)}_\lambda$ for the representation. 

Given this construction of $\varepsilon_\lambda$, it is straightforward to
compute its Chern
character
$\ch^*(\varepsilon_\lambda) \in HC^*[A(S^{N-1}_q)]$: It is the pull  back of
the Chern character
of the canonical generator of $K^0(\co)$.  An element of the cyclic
cohomology
$HC^0$ is a  trace. The degree $0$ part of the
Chern character of the canonical  generator of $K^0(\co)$ is given by
the identity map
$\co\to\co$,  which is trivially a trace. Pulling this back we find
$\ch^0(\varepsilon_\lambda) = \psi^{(1)}_\lambda : A(S^{N-1}_q)
\to \co$ which is also  a trace because it is a homomorphism to a
commutative algebra. These are distinct elements of $HC^0[A(S^{N-1}_q)]$ for different values of $\lambda$. However, because the Fredholm modules $\varepsilon_\lambda$ all lie in the same $K$-homology class, their Chern characters are all equivalent in periodic cyclic cohomology. Indeed, applying the periodicity operator once, the cohomology classes $S(\psi^{(1)}_\lambda) \in HC^2[A(S^{N-1}_q)]$ are all the same. For the computation of the pairing between $K$-theory and $K$-homology, any trace determining the same periodic cyclic cohomology class can be used. The most
symmetric choice of trace is given by  averaging $\psi^{(1)}_\lambda$
over $\lambda\in
S^1\subset \co$:
\[
\tau^0(a) := \oint_{S^1} \psi^{(1)}_\lambda(a) \frac{d\lambda}{2\pi i
\lambda} .
\]
The result is normalized, $\tau^0(1)=1$, and vanishes
on all the generators. The higher degree parts of $\ch^*(\varepsilon_\lambda)$
depend only on the $K$-homology class $[\varepsilon_\lambda]$ and can be
constructed from $\tau^0$ by the periodicity
operator.

\subsubsection{$K$-homology Generators for Even Spheres}
We will now construct an
element $[\mu_{\mathrm{ev}}] \in K^0[\C(S^{2n}_q)]$ by giving a
suitable even
Fredholm module
$\mu := (\Hi,F,\gamma)$.

Identify the Hilbert spaces for
the representations $\psi^{(2n)}_\pm$ given in \eqref{rep.even} by
identifying their bases, and
call this $\Hi$. The representation for
the Fredholm module
is
\[
\psi:=\psi^{(2n)}_+\oplus\psi^{(2n)}_-
\] acting on
$\Hi\oplus\Hi$. The grading operator and the Fredholm
operator are respectively,
\[
\gamma = \begin{pmatrix} 1&0
\\0&-1
\end{pmatrix}
,\quad
  F = \begin{pmatrix} 0&1\\1&0
\end{pmatrix}.
\]
It is obvious that $F$ is odd (since it anticommutes with $\gamma$)
and Fredholm (since
it is invertible).
The remaining property to check is that for any $a\in A(S^{2n}_q)$,
the commutator $[F,\psi(a)]_-$ is compact. Indeed,
\[
[F,\psi(a)]_- =
\begin{pmatrix} 0 & - \psi^{(2n)}_+(a) + \psi^{(2n)}_-(a) \\
\psi^{(2n)}_+(a) - \psi^{(2n)}_-(a) & 0
\end{pmatrix} .
\]
However, $\psi^{(2n)}_+(a) - \psi^{(2n)}_-(a) =
\psi^{(2n)}_+[a -
\sigma(a)]$ and $a-\sigma(a)$ is always proportional to a power of
$x_0$. Thus   this
is not only compact, it is trace class. This also shows that we have
(at least) a $1$-summable Fredholm module. This is in contrast to the fact that the analogous element of $K_0(S^{2n})$ for the undeformed sphere is given by a $2n$-summable Fredholm module.
\\
The Chern character \cite{Con94} $\ch^*(\mu_{\mathrm{ev}})$ has  a
component
in degree $0$, $\ch^{0}(\mu_{\mathrm{ev}})\in HC^0[A(S^{2n}_q)]$. The
element
$\ch^0(\mu_{\mathrm{ev}})$ is the trace
\begin{equation}\label{tau1}
\tau^1(a) := \Tr(\gamma \psi(a) ) = \Tr\left[\psi^{(2n)}_+(a) -
\psi^{(2n)}_-(a)\right].
\end{equation}
The higher degree parts of $\ch^*(\mu_{\mathrm{ev}})$
can be obtained via the periodicity operator.

For the sphere $S^{2}_q$ our Fredholm module coincides with the one 
constructed in
\cite{MasNakWat91}.

\subsubsection{$K$-homology Generators for Odd Spheres}
\label{se:khodd}
The element $[\mu_{\mathrm{odd}}] \in
K^1[\C(S^{2n+1}_q)]$ is most
easily given by an unbounded Fredholm
module.

Let the representation $\psi$ be
the direct integral (over $\lambda\in S^1$) of the
representations
$\psi^{(2n+1)}_\lambda$ given in \eqref{rep.odd}. The operator is the
unbounded operator
$D:= \lambda^{-1}\frac{d}{d\lambda}$.
\\
Referring to \eqref{rep.odd}, we see that the representative of $x_1$
is proportional to $\lambda$ and as a consequence,
\begin{subequations}
\begin{equation}
[D,\psi(x_1)]_- = \psi(x_1)
\end{equation}
whereas for $i>1$ the representative of $x_i$ does not involve
$\lambda$ and therefore
\begin{equation}
[D,\psi(x_i)]_- = 0, \quad i>0 .
\end{equation}
\end{subequations}
Since $a\mapsto [D,\psi(a)]_-$ is a derivation, this shows that
$[D,\psi(a)]_-$ is bounded for any $a \in A(S^{2n+1}_q)$. Note
however that for $n>0$ (i.~e., except for $S^1$) all eigenvalues of
$D$ have infinite degeneracy and therefore $D$ does not have compact
resolvent.

This triple can be converted in to a bounded  Fredholm module by
applying a cutoff
function to $D$. A convenient  choice is $F=\chi(D)$ where
\[
\chi(m) :=
\begin{cases} 1 &: m>0\\ -1  &: m\leq 0 .
\end{cases}
\] To be more explicit, use a Fourier  series basis for the Hilbert
space,
\[
\ket{k_0,k_1,\dots,k_n} :=
\lambda^{k_0}\ket{k_1,\dots,k_n} ,
\]
in which the  representation is given  by,
\begin{align*}
\psi(x_1)
\ket{k_0,\dots,k_n} &=
q^{-(k_1+\dots+k_n)}
\ket{k_0+1,\dots,k_n}
,  \\
\psi(x_1^*)
\ket{k_0,\dots,k_n} &=
q^{-(k_1+\dots+k_n)}
\ket{k_0-1,\dots,k_n},  \\
\psi(x_i)
\ket{k_0,\dots,k_n} &=  (1-q^{-2k_{i-1}})^{1/2} q^{-(k_i+\dots+k_n)}
\ket{k_0,\dots,k_{i-1}-1,\dots,k_n}, \\
\psi(x_i^*) \ket{k_0,\dots,k_n} &= (1-q^{-2(k_{i-1}+1)})^{1/2}
q^{-(k_i+\dots+k_n)}
\ket{k_0,\dots,k_{i-1}+1,\dots,k_n},
\end{align*}
for $i = 1, \dots, n$. The Fredholm operator is then given by
\[
F\ket{k_0,\dots,k_n} =
\chi(k_0) \ket{k_0,\dots,k_n}
.\]
The only condition to  check is that the commutator $[F,\psi(a)]_-$
is compact for any
$a\in  \C(S^{2n+1}_q)$. Since $a \mapsto [F,\psi(a)]_-$ is a
derivation, it  is
sufficient to check this on generators. One  finds
$[F,\psi(x_i)]_-=0$ for $i>1$  and
\begin{equation}
\label{F.com} [F,\psi(x_1)]_-
\ket{k_0,\dots,k_n} =
\begin{cases} 2q^{-(k_1+\dots+k_n)}
\ket{1,k_1,\dots,k_n} &:
k_0=0\\ 0 &: k_0\neq
0,
\end{cases}
\end{equation} which is indeed compact, and in
fact trace class.
\\
Thus, this is a $1$-summable Fredholm module. Again this is in contrast to the fact that the analogous element of
$K_1(S^{2n+1})$ for the undeformed sphere  is given by a $(2n+1)$-summable Fredholm module.
\\
Its Chern character \cite{Con94} begins with
$\ch^{\frac12}(\mu_{\mathrm{odd}})\in HC^1[A(S^{2n+1}_q)]$ which is
given by the cyclic
$1$-cocycle $\varphi$ defined
by
\begin{equation}
\label{phi.def}
\varphi(a,b) := \tfrac12
\Tr\left(\psi(a)[F,\psi(b)]_-\right)
.\end{equation}

The higher degree parts of $\ch^*(\mu_{\mathrm{odd}})$
can be obtained via the periodicity operator.

For the sphere $S^3_q \cong \SU_q(2)$ our Fredholm module coincides with the
one constructed in
\cite{MasNakWat90}.

\subsection{$K$-theory for Even  Spheres}
For $S^{2n}_q$ we construct
two classes in the $K$-theory group
$K_0[\C(S^{2n}_q)] \cong \Z^2$.

The first class is trivial. The element $[1]\in
K_0[\C(S^{2n}_q)]$ is
the equivalence class of $1\in \C(S^{2n}_q)$
which is of course an
idempotent. In order to compute the pairing
with $K$-homology, we need
the degree
$0$ part of its Chern
character,
$\ch_0[1]$, which is represented by the cyclic cycle $1$.

The second, nontrivial, class was presented in \cite{LanMad01}
\footnote{Again, for the sphere $S^{2}_q$ the idempotent \eqref{proj} 
was already presented
in \cite{MasNakWat91}.}
It is  given by an idempotent
$e_{(2n)}$ constructed from the unipotent \eqref{uni.even} as
\begin{equation}
\label{proj} e_{(2n)} =
\tfrac12(1+u_{(2n)}) .
\end{equation}
Its degree $0$ Chern character, $\ch_{0}(e_{(2n)})\in
HC_0[A(S^{2n}_q)]$, is
\begin{align*}
\ch_0(e_{(2n)}) = \tr(e_{(2n)}) & = 2^{n-1} + \tfrac12 \tr (u_{(2n)})
\\ & = 2^{n-1} + \tfrac12(q^{-1}-1)^n x_0 ,
\end{align*}
since the recursive definition \eqref{uni.even} of the unipotent
$u_{(2n)}$ shows that,
\[
\tr(u_{(2n)}) = (q^{-1}-1) \tr(u_{(2n-2)}) = (q^{-1}-1)^n x_0 .
\]

Now, we can pair these classes with the two $K$-homology elements
which we
constructed in Section~\ref{se:khom}. First,
\[
\langle \varepsilon_\lambda, [1]\rangle := \tau^0(1) =
1,
\]
which is hardly surprising. Second, the ``rank'' of  the idempotent
$e_{(2n)}$ is
\[
\langle \varepsilon_\lambda, e_{(2n)} \rangle
:= \tau^0(\ch_0 e_{(2n)}) = 2^{n-1}
.\]
Also not surprising is the ``degree'' of $[1]$,
\[
\langle \mu_{\mathrm{ev}}, [1]\rangle := \tau^1(1) =
\Tr\bigl[\psi^{(2n)}_+(1)-\psi^{(2n)}_-(1)\bigr] =
\Tr(1-1) = 0
.\]
The more complicated pairing is,
\begin{align*}
\langle \mu_{\mathrm{ev}}, e_{(2n)}\rangle  &:= \tau^1(\ch_0
e_{(2n)}) \\ &\:=
\Tr\circ\psi^{(2n)}_+\circ(1-\sigma)\left(2^{n-1} + \tfrac12
[q^{-1}-1]^nx_0\right) \\ &\:= (q^{-1}-1)^n \Tr[\psi^{(2n)}_+(x_0)]
.\end{align*} So, we need to compute
\[
\Tr[\psi^{(2n)}_+(x_0)] :=
\sum_{k_0=0}^\infty \dots \sum_{k_{n-1}=0}^\infty
q^{-(k_0+\dots+k_{n-1})} =
\left(\sum_{k=0}^\infty q^{-k}\right)^n = (1- q^{-1})^{-n}
.\]
And, in turn, we get,
\[
\langle \mu_{\mathrm{ev}}, e_{(2n)}\rangle = (-1)^n
.\]

The fact that the matrix of pairings,
\[
\begin{array}{c|cc}
& [1] & [e_{(2n)}] \\
\hline
[\varepsilon_\lambda] & 1 & 0 \\
{}[\mu_{\mathrm{ev}}] & 0 & (-1)^n
\end{array}
\]
is invertible over the integers proves that the classes
$[1],[e_{(2n)}] \in K_0[\C(S^{2n}_q)] \cong \Z^2$ and
$[\varepsilon_\lambda],[\mu_{\mathrm{ev}}]\in K^0[\C(S^{2n}_q)] \cong \Z^2$ are
generators of these groups.

Classically, the ``degree'' of the left-handed spinor bundle is $-1$.
So, the $K$-hom\-o\-logy class which correctly generalizes the classical
$K$-orientation class $[\mu]\in K_0(S^{2n})$ is actually
$(-1)^{n+1}[\mu_{\mathrm{ev}}]$.

\subsection{$K$-theory for Odd Spheres}
Again, define
$[1]\in
K_0[\C(S^{2n+1}_q)]$ as the equivalence class of
$1\in
\C(S^{2n+1}_q)$. The pairing with our element $[\varepsilon_\lambda]
\in
K^0[\C(S^{2n+1}_q)]$ is again,
\[
\langle \varepsilon_\lambda, [1]\rangle
:=\tau^0(1) = 1
.\]
There is no other independent
generator in $K_0[\C(S^{2n+1}_q)] \cong \Z$.

Instead, $K_1[\C(S^{2n+1}_q)] \cong \Z$ is nonzero. So we need to
construct a
generator there. An odd $K$-theory element is an equivalence class of
unitary
matrices over the algebra. We can construct an appropriate sequence
of unitary matrices recursively, just as we constructed the
unipotents and
idempotents.

Let
\begin{equation} \label{unit}
V_{(2n+1)} =
\begin{pmatrix} x_{n+1}  & q^{-1}
V_{(2n-1)} \\ -V_{(2n-1)}^* &
x_{n+1}^*
\end{pmatrix}
,
\end{equation}
with
$V_{(1)} = x_1$. By using the defining relations
\eqref{core} one
directly proves that this is unitary: $V_{(2n+1)}
V_{(2n+1)}^* =
V_{(2n+1)}^* V_{(2n+1)} = 1$.

In order to pair our
$K$-homology element $[\mu_{\mathrm{odd}}]\in
K^1[\C(S^{2n+1}_q)]$ with the unitary $V_{(2n+1)}$, we
need the lower degree part $\ch_{\frac12}(V_{(2n+1)}) \in
HC_1[A(S^{2n+1}_q)]$ of its Chern character. It is
given by the cyclic cycle,
\begin{align}
\ch_{\frac12}(V_{(2n+1)}) &:=
\tfrac12\tr\left(V_{(2n+1)} \otimes V_{(2n+1)}^* -
V_{(2n+1)}^*
\otimes V_{(2n+1)} \right) \nn \\ & \:=
\tfrac12(q^{-2}-1)^n  (x_1
\otimes x_1^* - x_1^*\otimes x_1) .
\end{align}
Now, compute the pairing,
\begin{align*}
\langle \mu_{\mathrm{odd}} , V_{(2n+1)}\rangle  &:= \langle
\varphi,
\ch_{\frac12}(V_{(2n+1)})\rangle \\ &\:= - (q^{-2}-1)^n
\varphi(x_1^*,x_1) \\ &\:= -
\tfrac12(q^{-2}-1)^n
\Tr\left(\psi(x_1^*)[F,\psi(x_1)]_-\right) .\end{align*}
We have already computed
$[F,\psi(x_1)]_-$ in eq.~\eqref{F.com}. From  that, we get
\[
\psi(x_1^*)[F,\psi(x_1)]_- \ket{k_0,\dots,k_n}  = \begin{cases} 2
q^{-2(k_1+\dots+k_n)} \ket{0,k_1,\dots,k_n} &: k_0=0\\ 0 &:
k_0\neq 0 .
\end{cases}
\]
Thus,
\begin{align*}
\Tr\left(\psi(x_1^*)[F,\psi(x_1)]_-\right)  &= \sum_{k_1=0}^\infty
\dots
\sum_{k_n=0}^\infty 2 q^{-2(k_1+\dots+ k_n)} \\ &= 2
\left(\sum_{k=0}^\infty q^{-2k}\right)^n \\ &= 2 (1- q^{-2})^{-n}
.\end{align*}
Finally, this gives
\[
\langle \mu_{\mathrm{odd}} , V_{(2n+1)}\rangle = (-1)^{n+1}
.\]
This proves the classes $[V_{(2n+1)}]\in K_1[\C(S^{2n+1}_q)]$ and
$[\mu_{\mathrm{odd}}]\in K^1[\C(S^{2n+1}_q)]$ are nonzero and that
neither may be a multiple of another class. Thus $[V_{(2n+1)}]$ and
$[\mu_{\mathrm{odd}}]$ are indeed generators of these groups.

\section{Final remarks}
Recently, in \cite{ChaPal02}, a $3$-summable spectral triple was
constructed for $\SU_q(2)$; this has been thoroughly analyzed in 
\cite{Con02}
in the context of the noncommutative local index formula of 
\cite{ConMos95}. Also,
a $2$-summable spectral triple on $\SU_q(2)$ was constructed in
\cite{ChaPal02b} together with a spectral triple on the  spheres $S^2_{qc}$
of Podle{\`s} \cite{Pod87}.

It would be interesting to compare our $1$-summable Fredholm modules
for the sphere $S^3_q \cong \SU_q(2)$ and for the equator sphere $S^2_q \cong
S^2_{q\infty}$ with these spectral triples. It seems likely that for the same
spheres  they determine the same $K$-homology class.

It would also be interesting, although much more challenging, to extend (some
of) the analysis of
\cite{{Con02}} to all spheres, notably odd ones, $S^{2n+1}_q$, with a suitable
modification of the unbounded Fredholm modules constructed in
Section~\ref{se:khodd}.

\section{Acknowledgments}
We are grateful to J.~Cuntz and N.~Higson for directing us to some
references and to F.~Bonechi,  G.~Fiore, L.~Dabrowski, J.~Madore, and D.~Perrot for useful  conversations. J.~Varilly's remarks on the previous version have led to significant improvements.

\providecommand{\href}[2]{#2}\begingroup\raggedright
\endgroup

\end{document}